\documentclass{elsart}
\usepackage{amssymb}
\usepackage{amsfonts}
\usepackage{graphics}
\usepackage{graphicx}
\begin{document}

\begin{frontmatter}



\title{Maximum edges possible in a graph for restricted independence number, maximum degree, and maximum matching size}
{\author{Niraj Khare},
\ead{nirajkhare@math.ohio-state.edu}
\author{Nishali Mehta},
\ead{nishali@math.ohio-state.edu}
\author{Naushad Puliyambalath\corauthref{cor}}
\corauth[cor]{Corresponding author}}\footnote{Authors share equal credit for the work.}
\ead{pasha@math.ohio-state.edu}
\address{The Ohio State University, Columbus, Ohio, USA }
\begin{abstract}
This article provides sharp bounds for the maximum number of edges possible in a simple graph with restricted values of two of the three parameters, namely, maximum matching size, independence number and maximum degree. We also construct extremal graphs that achieve the edge bounds in all cases. We further establish uniqueness of these extremal graphs whenever they are unique. 
\end{abstract}
\end{frontmatter}
\section{Introduction}
By a graph, we shall mean a simple graph, i.e., a graph with no loop and no multiple edges. We first fix some notation. For a graph $G$, $E(G)$ and $V(G)$ would denote the edge set and the vertex set of $G$, respectively. $\alpha(G)$, $\Delta(G)$, and $\nu(G)$ would denote the independence number of $G$, the maximum degree of any vertex in $G$, and the size of a maximum matching in $G$, respectively. For $x \in V(G)$, $deg_G(x)$ would denote the degree of the vertex $x$ and $G\setminus x$ would denote the induced subgraph on $V(G)\setminus \{x\}$. $\chi(G)$ and $\chi'(G)$ would denote the chromatic number and the chromatic index of $G$, respectively.\\
We now consider the problem of finding for a graph $G$ the maximum of $|E(G)|$ when two of the three parameters $\alpha(G)$, $\Delta(G)$ and $\nu(G)$ are known. We will establish shortly that the problem is well-founded as crude bounds exist.\newline
The problem of finding a precise upper bound for the maximum number of edges in a graph $G$, where $\nu(G)$, $\Delta(G)$ and $|V(G)|$ are known, is discussed in \cite{AH}, \cite{CH}. A simpler proof is provided in \cite{BK} for a precise bound on $|E(G)|$ when only $\nu(G)$ and $\Delta(G)$ are known. The works \cite{AH}, \cite{CH} and  \cite{BK} were inspired by well-known work of Erd\H{o}s-Rado \cite{ER} and consequently do not consider $\alpha(G)$ as one of the parameters. This article finds the precise upper bound on $|E(G)|$ when two of the three parameters: $\alpha(G)$, $\Delta(G)$ and $\nu(G)$ are known. It further finds the complete set of conditions so that the graph achieving the upper bound on $|E(G)|$ is unique. 
\section{Results}
Our main aim in this article is to prove the following results.
\begin{defn}
For $\alpha, \nu\in \mathbb{Z}^+$, we define two graphs $G_{\alpha,\nu}$ and $H_{\alpha,\nu}$ as follows:
(a) $G_{\alpha, \nu}$ contains $\alpha$ connected components: $K_{2\nu+1}$ and $\alpha-1$ isolated vertices.\\
Note: $|E(G_{\alpha, \nu})|={2\nu+1 \choose 2}=\nu(2\nu+1)$, $\alpha(G_{\alpha, \nu})=\alpha$, and $\nu(G_{\alpha, \nu})=\nu$.\\
(b) $H_{\alpha, \nu}$: Start with the complete bipartite graph $K_{\alpha, \nu}$. Add ${\nu \choose 2}$ edges to the part with $\nu$ vertices.\\
Note: $|E(H_{\alpha, \nu})|=\alpha\nu+{\nu \choose 2}$, $\alpha(H_{\alpha, \nu})=\alpha$, and, when $\alpha \geq \nu$, $\nu(H_{\alpha, \nu})=\nu$.
\end{defn}
\begin{thm}\label{alpha-nu}
For all $\alpha, \nu \in \mathbb{Z}^+$, let $G$ be a graph with $\alpha(G)\leq \alpha$ and $\nu(G)\leq \nu$.\\
(a) If $2\alpha < 3(\nu+1)$, then $|E(G)|\leq {2\nu+1\choose 2}$ with equality holding iff $G=G_{\alpha, \nu}$.\\
(b) If $2\alpha > 3(\nu+1)$, then $|E(G)|\leq \alpha\nu+{\nu\choose 2}$ with equality holding iff $G=H_{\alpha, \nu}$.\\
(c) If $2\alpha = 3(\nu+1)$, then $|E(G)|\leq {2\nu+1\choose 2}=\alpha\nu+{\nu\choose 2}$ with equality holding iff $G \in \{G_{\alpha, \nu},H_{\alpha, \nu}\}$.\\ 
\end{thm}

\begin{thm}\label{alpha-delta}
Let $G$ be a graph with $\Delta(G)\leq \Delta$ and $\alpha(G)\leq \alpha$. Then $|E(G)|\leq \alpha{\Delta+1\choose 2}$ and the bound is sharp. Furthermore, the unique graph that attains the edge bound consists of $\alpha$ components where each component is the complete graph $K_{\Delta+1}$.
\end{thm}

A shorter proof of the following theorem can be found in \cite{BK}:

\begin{thm}\label{delta-nu}
For $\Delta, \nu \in \mathbb{Z}^+$, let $G$ be a graph with $\Delta(G) = \Delta$ and $\nu(G) = \nu$, then
$$|E(G)|\leq \Delta\nu+\lfloor\frac{\nu}{\lceil\frac{\Delta}{2}\rceil}\rfloor \lfloor\frac{\Delta}{2}\rfloor $$ and the edge bound is sharp for all $\Delta$ and $\nu$.
\end{thm}

We look at those graphs that achieve the edge bound in Theorem \ref{delta-nu} and prove the following:

\begin{thm} \label{delta-nu2}
For $\Delta, \nu \in \mathbb{Z}^+$, let $G$ be a graph with $\Delta(G) = \Delta$ and $\nu(G) = \nu$. If $G$ is a graph with no isolated vertices and $|E(G)| = \Delta\nu+\lfloor\frac{\nu}{\lceil\frac{\Delta}{2}\rceil}\rfloor \lfloor\frac{\Delta}{2}\rfloor$, $G$ is unique if and only if $\lceil\frac{\Delta}{2}\rceil$ divides $\nu$ or $\nu=1$.
\end{thm}

\section{Preliminaries}
In this section we list well-known results which we will use in this paper. We first start with vertex coloring related results. These results can be found in any standard text book on graph theory. In particular readers can find them in \cite{LP} or \cite{DW}.

\begin{prop}\label{d+1 vertex-bound}
$\chi(G)\leq \Delta(G)+1$ for every graph $G$.
\end{prop}
\begin{thm} ({\bf{Brooks}}) \label{Brooks thm}
Let $G$ be a connected graph. $\chi(G)=\Delta(G)+1$ if and only if $G$ is either an odd cycle or a complete graph.
\end{thm}
Now we state matching and edge coloring related results. We define a factor-critical graph and state Gallai's Lemma, which is crucial to the following discussion. An elegant proof of Gallai's Lemma can be found in \cite{LP}.
\begin{defn}\label{factor-critical}
A connected graph $G$ is called factor-critical if and only if $G\setminus x$ has a perfect matching for all $x \in V(G)$.
\end{defn}
Note that if $G$ is a factor-critical graph, then $|V(G)|=2\nu(G)+1$.
\begin{lem}{\bf{(Gallai)}}\label{Gallai's lemma}
Let $G$ be a connected graph. If $\nu(G\setminus x)=\nu(G)$ for all $x \in V(G)$, then $G$ is a factor-critical graph.
\end{lem}
Finally, we recall the celebrated result of Vizing about chromatic index of a simple graph.
\begin{thm} ({\bf{Vizing, Gupta}})\label{vizing}
$\chi'(G)\leq \Delta(G)+1$ for every graph $G$.
\end{thm}
We borrow definitions from \cite{LP} to state the following two results. We shall use these definitions in section $4$ as well.
For any graph $G$, we define a set $S_G \subseteq V(G)$ by
\begin{displaymath}
S_G=\{v\in V(G)\,:\, \nu(G\backslash v) = \nu(G)-1\}.
\end{displaymath}
In other words, if $v \in S_G$, then $v$ must be covered by every maximum matching of $G$. 
\begin{eqnarray*}
D(G) &=& V(G)\setminus S_G, \\
A(G) &=& \{v \in S_G \quad| \quad v \textrm{ is adjacent to a vertex in } D(G) \}, \textrm{ and} \\
C(G) &=& S_G\setminus A(G).
\end{eqnarray*}
\begin{thm} ({\bf{Edmond-Gallai Structure Theorem}}) \label{Ed-Gal}
If $G$ is a simple graph and $A(G)$, $C(G)$ and $D(G)$ are defined as above, then:\\
(a) the components of the graph induced by $D(G)$ are factor critical,\\
(b) the subgraph induced by $C(G)$ has a perfect matching,\\
(c) any maximum matching of G contains a near perfect matching of $D(G)$, a perfect matching of $C(G)$ and matches all vertices of $A(G)$ with vertices in distinct components of $D(G)$.
\end{thm}
\begin{lem} ({\bf{Stability Lemma}}) \label{stability}
If $G$ is a simple graph and $A(G)$, $C(G)$ and $D(G)$ are defined as above.\\
(i) Let $u \in A(G)$. Then $A(G\setminus u)=A(G) \setminus \{u\}$, $C(G \setminus u)= C(G)$ and $D(G \setminus u)= D(G)$.\\
(ii) Let $u \in C(G)$. Then $A(G) \subseteq A(G\setminus u)$, $C(G \setminus u)\subseteq C(G)\setminus \{u\}$ and $D(G) \subseteq D(G \setminus u)$.\\
(iii) Let $u \in D(G)$. Then $A(G\setminus u)\subseteq A(G)$, $C(G)\subseteq C(G \setminus u)$ and $D(G \setminus u)\subseteq D(G)\setminus \{u\}$.
\end{lem}
We next find the edge bounds mentioned earlier by considering all pairs of two of the three parameters: independence number $\alpha$, maximum degree $\Delta$, and maximum matching size $\nu$. 
%****************************************************************************
%*****************************************************************************
\section{Graphs with restricted $\alpha$ and $\nu$}
For this section we fix $\alpha,\nu \in\mathbb{Z}^{+}$. Let $G$ be a graph with $\alpha(G)=\alpha$ and $\nu(G)= \nu$. For a matching $M$ in $G$, let $V(M)$ denote the vertices in $G$ that are incident to edges in $M$. Note that $|V(M)|=2|M|$. If $M$ is a maximal matching in G, the set of vertices $V(G)\backslash V(M)$ is independent, so $|V(G)|-|V(M)|$ can be at most $\alpha$. Such a graph $G$ must satisfy $|V(G)| \le 2|M|+\alpha\le 2\nu+\alpha$ and $|E(G)| \le {2\nu+\alpha \choose 2}$. Define a function $e_{1}: \mathbb{Z}^{+} \times \mathbb{Z}^{+} \to \mathbb{Z}^{+}$ and a set of graphs $E_{1}(\alpha, \nu)$ as follows:
\begin{displaymath}
e_{1} (\alpha, \nu) = \max \{ |E(G)| : \alpha(G)=\alpha, \nu(G)=\nu \} \le {2\nu+\alpha \choose 2},
\end{displaymath}
\begin{displaymath}
E_{1}(\alpha, \nu) = \{ G : \alpha(G)=\alpha, \nu(G)=\nu, |E(G)| = e_{1}(\alpha, \nu) \}.
\end{displaymath}
We define two graphs $G_{\alpha,\nu}$ and $H_{\alpha,\nu}$ as follows:
\begin{itemize}
\item $G_{\alpha, \nu}$ contains $\alpha$ connected components: $K_{2\nu+1}$ and $\alpha-1$ isolated vertices. It is clear that $|E(G_{\alpha, \nu})|={2\nu+1 \choose 2}$, $\alpha(G_{\alpha, \nu})=\alpha$, and $\nu(G_{\alpha, \nu})=\nu$.

\item $H_{\alpha, \nu}$: Start with the complete bipartite graph $K_{\alpha, \nu}$. Add ${\nu \choose 2}$ edges to the part with $\nu$ vertices.
It is clear that $|E(H_{\alpha, \nu})|=\alpha\nu+{\nu \choose 2}$, $\alpha(H_{\alpha, \nu})=\alpha$, and, when $\alpha \geq \nu$, $\nu(H_{\alpha, \nu})=\nu$.
\end{itemize}

Since $G_{\alpha, \nu}$ and $H_{\alpha, \nu}$ have the appropriate $\alpha$ and $\nu$ values, we have:

\begin{itemize}
\item for $\alpha < \nu$,
\begin{equation}\label{eqn1}
e_1(\alpha, \nu) \ge \nu(2\nu+1),
\end{equation}
\item for $\alpha \geq \nu$,

\begin{equation}\label{eqn2}
e_1(\alpha, \nu) \ge \max\left\{\nu(2\nu+1),\; \nu\alpha+{\nu \choose 2}\right\}.
\end{equation}
\end{itemize}

We will show that equality holds in both (\ref{eqn1}) and (\ref{eqn2}) and that $G_{\alpha,\nu}$ and $H_{\alpha,\nu}$ are the only extremal graphs.

\begin{prop}\label{unsat}
If $G \in E_{1}(\alpha, \nu)$ and $\alpha>1$, then $|V(G)| \geq 2\nu+2$.
\end{prop}
{\sl Proof: } Suppose $G \in E_{1}(\alpha, \nu)$, $\alpha>1$, and $|V(G)| \leq 2\nu+1$. $G$ has at most ${2\nu+1 \choose 2}$ edges. Also $|E(G)|=e_{1}(\alpha, \nu) \geq {2\nu+1 \choose 2}$, so $|E(G)| = {2\nu+1 \choose 2}$. Thus $G=K_{2\nu+1}$, implying $\alpha=1$, which is a contradiction.\qed

\begin{prop}\label{every}
If $G\in E_{1}(\alpha, \nu)$ and $v\in S_G$, then $v$ is adjacent to every other vertex of $G$.
\end{prop}
{\sl Proof: } Let $G\in E_{1}(\alpha, \nu)$, and $v\in S_G$. Suppose $u \in V(G)$, $u \neq v$, and $(u,v) \notin E(G)$. Define a new graph $G'$ as follows: $V(G')=V(G)$ and $E(G')=E(G) \cup \{(u,v)\}$. If $\alpha(G')=\alpha-1$, modify $G'$ by adding a new vertex to $V(G')$ with no additional edges.\\
By the Edmond-Gallai Structure Theorem (Theorem \ref{Ed-Gal}) $\nu(G')=\nu(G)$ as $v\in S_G=A(G)\cup C(G)$. So we get $\nu(G')=\nu(G)$, $\alpha(G')=\alpha$ and $|E(G')|=|E(G)|+1$. But this contradicts the fact that $|E(G)|=e_{1}(\alpha, \nu)$. Therefore, $v$ must be adjacent to every other vertex of $G$. \qed 
\begin{prop}\label{nured}
If $G\in E_{1}(\alpha, \nu)$, then $\nu(G\backslash S_G)=\nu-|S_G|$. 
\end{prop}
{\sl Proof: } Since $G\in E_{1}(\alpha, \nu)$, we know $|V(G)|\geq 2\nu+1$. Therefore, $D(G)$ is not an empty set. So by previous proposition $S_G=A(G)$. By the Stability Lemma (Lemma \ref{stability}), we have $\nu(G\backslash A(G))=\nu-|A(G)|=\nu-|S_G|$.\qed
\begin{thm} \label{alphanu}
For all $\alpha, \nu \in \mathbb{Z}^{+}$:
\begin{itemize}
\item[(a)] If $2\alpha < 3(\nu+1)$, then $e_{1}(\alpha, \nu)={2\nu+1\choose 2}$ and $E_{1}(\alpha, \nu)=\{G_{\alpha, \nu}\}$.
\item[(b)] If $2\alpha > 3(\nu+1)$, then $e_{1}(\alpha, \nu)=\alpha\nu+{\nu\choose 2}$ and $E_{1}(\alpha, \nu)=\{H_{\alpha, \nu}\}$.
\item[(c)] If $2\alpha = 3(\nu+1)$, then $e_{1}(\alpha, \nu)={2\nu+1\choose 2}=\alpha\nu+{\nu\choose 2}$ and\\ $E_{1}(\alpha, \nu)=\{G_{\alpha, \nu},H_{\alpha, \nu}\}$. 
\end{itemize}
\end{thm}
{\sl Proof: } Given $\alpha, \nu \in \mathbb{Z}^{+}$, let $G \in E_{1}(\alpha, \nu)$.

If $\alpha=1$, then for any $\nu \geq 1$, $2\alpha=2 < 3(\nu+1)$. Also, $G$ must be the complete graph $K_{2\nu+1}$, which is $G_{1,\nu}$. Then $e_{1}(1, \nu)=\nu(2\nu+1)$ and $E_{1}(1, \nu)=\{G_{1, \nu}\}$.

Assume $\alpha>1$. Let $s=|S_G|$, and $G'= G\backslash S_G$. Let $G_1, G_2, \dots, G_k$ be the connected components of $G'$, and $\nu_i=\nu(G_i)$ for $1 \leq i \leq k$. Note that $k \leq \alpha$ and by the Edmond-Gallai Structure Theorem (Theorem \ref{Ed-Gal}) each $G_i$ is factor-critical. By Proposition \ref{nured}, we have $\nu_1+\nu_2+\cdots+\nu_k=\nu-s$ and $0\le s \le \nu$. By Gallai's Lemma (Lemma \ref{Gallai's lemma}), $G_i$ has exactly $2\nu_i+1$ vertices. Thus, $G_i$ has at most $\nu_i(2\nu_i+1)$ edges. There are three types of edges in $G$, namely $(x,y)$ where (i) $x,y \in S_G$, (ii) $x \in S_G$, $y \notin S_G$, and (iii) $x,y \notin S_G$.  Applying Proposition \ref{every}, we estimate $|E(G)|$ as follows:
\begin{eqnarray*}
|E(G)| &\le& {s\choose 2}+s|V(G')|+|E(G')|\\
&=& {s\choose 2}+s\left(\sum_{i=1}^k|V(G_i)|\right)+\sum_{i=1}^k |E(G_i)|\\
&\le& {s\choose 2}+s\left(\sum_{i=1}^k(2\nu_i+1)\right)+\sum_{i=1}^k\nu_i(2\nu_i+1)\\
&=& {s\choose 2}+s\left(\sum_{i=1}^k(2\nu_i+1)\right)+\left( 2\sum_{i=1}^k\nu_i^2 + \sum_{i=1}^k\nu_i\right) \\
&\le& {s\choose 2}+s\left(2\sum_{i=1}^k\nu_i+\sum_{i=1}^k 1\right)+\left[ 2\left( \sum_{i=1}^k\nu_i\right)^2 + \sum_{i=1}^k\nu_i \right] \\
&=& \frac{1}{2}s(s-1)+s(2(\nu - s)+k)+[2(\nu-s)^2+(\nu-s)]\\
&\le& \frac{1}{2}s(s-1)+s(2(\nu - s)+\alpha)+[2(\nu-s)^2+(\nu-s)]\\
&=& (\frac{1}{2}s^2-\frac{1}{2}s)+(2\nu s-2s^2+\alpha s) + (2\nu^2-4\nu s +2s^2+\nu- s) \\
&=& \frac{1}{2}s^2-\frac{3}{2}s + \alpha s + 2\nu^2-2\nu s+\nu- s.
\end{eqnarray*}
Rearranging terms, we get:
\begin{eqnarray} \label{smallalpha}
|E(G)|&\le& 2\nu^2+\nu - \frac{1}{2}\nu s + \frac{1}{2}s^2 - \frac{3}{2} \nu s - \frac{3}{2}s + \alpha s \nonumber \\
&=&\nu(2\nu+1)-\frac{1}{2}s(\nu-s)-\frac{1}{2}s\left(3(\nu+1)-2\alpha\right),
\end{eqnarray}
and
\begin{eqnarray} \label{bigalpha}
|E(G)| &\le& \frac{1}{2}s^2 - \frac{3}{2}s + \alpha s + \frac{1}{2}\nu^2 + \frac{3}{2}\nu^2 - \frac{1}{2} \nu s - \frac{3}{2}\nu s - \frac{1}{2}\nu + \frac{3}{2}\nu + \alpha \nu - \alpha \nu \nonumber \\
&=& \alpha \nu  + \frac{1}{2}(\nu^2 - \nu) - \frac{1}{2}(\nu s - s^2) - \frac{1}{2}(2\alpha\nu-3\nu^2-3\nu) + \frac{1}{2}(2\alpha s-3\nu s-3s) \nonumber \\
&=&\alpha\nu+{\nu\choose 2}-\frac{1}{2}s(\nu-s)-\frac{1}{2}(\nu-s)(2\alpha-3(\nu+1)).
\end{eqnarray}
Since $0 \le s \le \nu$, we have $\frac{1}{2}s(\nu-s) \ge 0$. Now we have three cases.
\begin{itemize}
\item[(a)] Suppose $2\alpha < 3(\nu+1)$. By (\ref{smallalpha}), $|E(G)| \le \nu(2\nu+1)$. Also, $|E(G)| =e_{1}(\alpha, \nu) \ge \nu(2\nu+1)$. Thus the above inequalities are all equalities. This happens if and only if $s=0$, $k=\alpha$, and there is $i \in \{ 1,2,\dots,k \}$ such that $\nu_i=\nu$, $|E(G_i)|=\nu_i(2\nu_i+1)$, and, for $j \neq i$, $\nu_j = 0$. This occurs when $G=G_{\alpha, \nu}$. Therefore $e_{1}(\alpha, \nu) = \nu(2\nu+1)$ and $E_{1}(\alpha, \nu)=\{G_{\alpha, \nu}\}$.

\item[(b)] Suppose $2\alpha > 3(\nu+1)$. Note that this implies $\alpha > \nu$. By (\ref{bigalpha}), $|E(G)| \le \alpha\nu+{\nu\choose 2}$. Also, $|E(G)| =e_{1}(\alpha, \nu) \ge \alpha\nu+{\nu\choose 2}$. Thus the above inequalities are all equalities. This happens if and only if $s=\nu$, $k=\alpha$, and $\nu_i=0$ for all $i$. This occurs when $G=H_{\alpha, \nu}$. Therefore $e_{1}(\alpha, \nu) = \alpha\nu+{\nu\choose 2}$ and $E_{1}(\alpha, \nu)=\{H_{\alpha, \nu}\}$.

\item[(c)] Suppose $2\alpha = 3(\nu+1)$. Note that this implies $\alpha > \nu$. In this case, we have $\nu(2\nu+1)=\alpha\nu+{\nu\choose 2}$. As in the previous two cases, the above inequalities are all equalities. So either $s=0$ or $s=\nu$.

If $s=0$, there is $i \in \{ 1,2,\dots,k \}$ such that $\nu_i=\nu$, $|E(G_i)|=\nu_i(2\nu_i+1)$, and, for $j \neq i$, $\nu_j = 0$. This occurs when $G=G_{\alpha, \nu}$.

If $s=\nu$, then, by Proposition \ref{nured}, $\nu(G')=0$. This means $V(G)\backslash S_G$ is an independent set in $G$. By Proposition \ref{every}, all vertices in $S_G$ are adjacent to all other vertices. This occurs when $G=H_{\alpha, \nu}$.

Therefore $e_{1}(\alpha, \nu) = {2\nu+1\choose 2}=\alpha\nu+{\nu\choose 2}$ and $E_{1}(\alpha, \nu)=\{G_{\alpha, \nu}, H_{\alpha, \nu}\}$.
\end{itemize}
This completes the proof.\qed

For given $\alpha, \nu  \in \mathbb{Z}^+$, let $G$ be a graph such that $\alpha(G) \le\alpha$ and $\nu(G)\le\nu$. Keeping one of $\alpha$ or $\nu$ fixed, $e_2(\alpha,\nu)=\max\left\{{2\nu+1\choose 2},\; \nu\alpha+{\nu \choose 2}\right\}$ is an increasing function in the other variable. Thus we have
\begin{displaymath}
|E(G)| \leq e_2(\alpha(G),\nu(G))\le e_2(\alpha,\nu)=\max\left\{{2\nu+1\choose 2},\; \nu\alpha+{\nu \choose 2}\right\}
\end{displaymath}
with equality holding if and only if $\alpha(G)=\alpha, \nu(G)=\nu$, and $G\in\{ G_{\alpha,\nu},H_{\alpha,\nu}\}$. This proves Theorem \ref{alpha-nu}.

We would like to mention here that Erd\H{o}s and Gallai considered the problem of finding the maximum number of edges in a graph with fixed vertex number $n$ and matching number $\nu$ in \cite{EG}. Their result is very similar to Theorem \ref{alpha-nu}.

\begin{thm}({\bf{Erd\H{o}s-Gallai}}) \label{Erdos-Gal}
For all $n, \nu \in \mathbb{Z}^+$, let $G$ be a graph with $|V(G)|\leq n$, $\nu(G)\leq \nu$, and $2\nu+2 \leq n$. Then
$$ |E(G)| \leq \max \left\{ {2\nu+1\choose 2}, (n-\nu)\nu+{\nu\choose 2} \right\}. $$
\end{thm}

A short proof of this result can be found in \cite{AF}.

In general, $\alpha \leq n - \nu$. When $\alpha = n - \nu$, the two bounds are identical. The fact that this is true is nontrivial. Using Theorem \ref{alpha-nu}, we can find a restriction on the graphs that achieve the maximum edge bound in Theorem \ref{Erdos-Gal}:

\begin{lem} \label{alpha-nu-n}
For all $n, \nu \in \mathbb{Z}^+$, let $G$ be a graph with $|V(G)|\leq n$, $\nu(G)\leq \nu$, and $2\nu+2 \leq n$. If $\alpha(G) < n - \nu$ and $ {2\nu+1\choose 2} < (n-\nu)\nu+{\nu\choose 2} $, then $$|E(G)| < (n-\nu)\nu+{\nu\choose 2}.$$
\end{lem}

{\sl Proof: }
By Theorem \ref{alpha-nu}, $|E(G)| \leq \max \left\{ {2\nu+1\choose 2}, \alpha\nu+{\nu\choose 2} \right\} < (n-\nu)\nu+{\nu\choose 2}.$ \qed

That is, a graph $G$ that achieves the maximum edge bound in Theorem \ref{Erdos-Gal} must satisfy $\alpha(G) = |V(G)| - \nu(G)$ or $\alpha(G)=|V(G)|-2\nu(G)$. Note that each graph $G \in E_1(\alpha, \nu)$ satisfies $\alpha(G) = |V(G)| - \nu(G)$ or $\alpha(G)=n-2\nu(G)$. Therefore, we can take this one step further: 

\begin{thm} \label{alpha-nu-n2}
For all $n, \nu \in \mathbb{Z}^+$, let $G$ be a graph with $|V(G)|\leq n$, $\nu(G)\leq \nu$, and $2\nu+2 \leq n$. If $|E(G)| = \max \left\{ {2\nu+1\choose 2}, (n-\nu)\nu+{\nu\choose 2} \right\}$, then $$G \in E_1(n-\nu, \nu) \cup E_1(n-2\nu, \nu).$$
\end{thm}

%*************************************************************************************
%************************************************************************************
\section{Graphs with restricted $\alpha$ and $\Delta$}
For this section, we fix $\alpha, \Delta\in \mathbb{Z}^{+}$. Let $G$ be a graph with $\alpha(G)=\alpha$ and $\Delta(G)= \Delta$. Consider the vertex chromatic number $\chi(G)$, the size of the smallest proper vertex coloring of $G$. There is a way to partition the vertices of $G$ into $\chi(G)$ color classes, where each color class is an independent set of $G$. By Proposition \ref{d+1 vertex-bound}, $\chi(G) \le \Delta+1$ and each color class has size at most $\alpha$. Therefore $|V(G)| \le \alpha(\Delta+1)$ and hence $|E(G)| \le {\alpha(\Delta+1)\choose 2}$. Define a function $e_{2}: \mathbb{Z}^{+} \times \mathbb{Z}^{+} \to \mathbb{Z}^{+}$ and a set of graphs $E_{2}(\alpha, \Delta)$ as follows:
\begin{displaymath}
e_{2} (\alpha, \Delta) = \max \{ |E(G)| : \alpha(G)=\alpha, \Delta(G)=\Delta \} \le {\alpha(\Delta+1)\choose 2},
\end{displaymath}
\begin{displaymath}
E_{2}(\alpha, \Delta) = \{ G : \alpha(G)=\alpha, \Delta(G)=\Delta, |E(G)| = e_{2} (\alpha, \Delta) \}.
\end{displaymath}
Define the graph $F_{\alpha, \Delta}$ to be the graph consisting of $\alpha$ connected components, where each component is the complete graph $K_{\Delta+1}$.

\begin{thm} \label{alphadelta}
For all $\alpha, \Delta \in \mathbb{Z}^{+}$, $e_{2} (\alpha, \Delta) = \alpha{\Delta+1\choose 2}$ and $E_{2}(\alpha, \Delta) = \{ F_{\alpha, \Delta} \}$.
\end{thm}
{\sl Proof. } The graph $F_{\alpha, \Delta}$ satisfies $\Delta(F_{\alpha, \Delta})=\Delta$ and $\alpha(F_{\alpha, \Delta})=\alpha$, implying 
\begin{displaymath}
e_{2} (\alpha, \Delta)\ge |E(F_{\alpha, \Delta})|=\alpha{\Delta+1\choose 2}.
\end{displaymath}
Also, if $G \in E_{2}(\alpha, \Delta)$, then, as observed above, $|V(G)| \le \alpha(\Delta+1)$ and
\begin{displaymath}
e_{2} (\alpha, \Delta)=|E(G)|=\sum_{v\in V(G)}\frac{1}{2}deg_{G}(v)\le \frac{\Delta}{2}|V(G)|\le\alpha{\Delta+1\choose 2}.
\end{displaymath}
Thus $e_{2} (\alpha, \Delta) = \alpha{\Delta+1\choose 2}$.\\
To prove uniqueness of $G$, let  $G \in E_{2}(\alpha, \Delta)$ and $G_1, G_2,\ldots, G_r$ be the connected components of $G$. Then 
\begin{eqnarray*}
\alpha\Delta(\Delta+1)&=&2|E(G)| =\sum_{i=1}^r2|E(G_i)|\\
&\le&\sum_{i=1}^r|V(G_i)|\Delta(G_i)\le\sum_{i=1}^r\chi(G_i)\alpha(G_i)\Delta(G_i) \\
&\le& \sum_{i=1}^r(\Delta(G_i)+1)\alpha(G_i)\Delta(G_i)\le \sum_{i=1}^r(\Delta+1)\alpha(G_i)\Delta \\
&=&\alpha\Delta(\Delta+1).
\end{eqnarray*}
All the above inequalities must be equalities and for each $G_i$ we have

\begin{tabular}{rl}
(i) & $\Delta(G_i)=\Delta,$ \\
(ii) & $\chi(G_i)=\Delta(G_i)+1=\Delta+1$, and \\
(iii) & $|V(G_i)|=\chi(G_i)\alpha(G_i)=(\Delta+1)\alpha(G_i)$.
\end{tabular}

By Brooks' Theorem (Theorem \ref{Brooks thm}), condition (ii) implies that every $G_i$ is either an odd cycle or the complete graph $K_{\Delta+1}$. Condition (iii) implies that if $G_i$ is an odd cycle then $\Delta=2$ and $G_i=K_3=K_{\Delta+1}$ (for if $G_i=C_{2n+1}$, then $\Delta=\Delta(G_i)=2$, $|V(G_i)|=2n+1$ and $\alpha(G_i)=n$ implying $2n+1=3n$). Thus every $G_i$ is $K_{\Delta+1}$ and there must be $\alpha$ of them, proving that $G$ is $F_{\alpha,\Delta}$.\qed

For given $\alpha, \Delta \in \mathbb{Z}^{+}$, let $G$ be a graph such that $\Delta(G)\le\Delta$ and $\alpha(G) \le\alpha$. Keeping one of $\alpha$ or $\Delta$ fixed, $e_1(\alpha,\Delta)=\alpha{\Delta+1\choose 2}$ is an increasing function in the other variable. Thus we have
\begin{displaymath}
|E(G)| \leq e_1(\alpha(G),\Delta(G))\le e_1(\alpha,\Delta)=\alpha{\Delta+1\choose 2}
\end{displaymath}
with equality holding if and only if $\alpha(G)=\alpha, \Delta(G)=\Delta$, and $G=F_{\alpha,\Delta}$. This proves Theorem \ref{alpha-delta}.
%****************************************************************************
%*****************************************************************************
\section{Graphs with restricted $\Delta$ and $\nu$}
For this section we fix $\Delta,\nu \in\mathbb{Z}^{+}$. We consider graphs with \underline{no} \underline{isolated} \underline{vertices} only. Let $G$ be a graph with $\Delta(G) \le \Delta$ and $\nu(G)\le \nu$. For every proper coloring of the edges of $G$, each color class is a matching and has at most $\nu$ edges. Thus using a minimum proper edge coloring, we get $|E(G)|\leq \chi'(G)\nu(G)$. By Theorem \ref{vizing}, $\chi'(G)\leq \Delta(G)+1$. Thus $G$ has at most $(\Delta+1)\nu$ edges.  

Define a function $e_{3}: \mathbb{Z}^{+} \times \mathbb{Z}^{+} \to \mathbb{Z}^{+}$ and a set of graphs $E_{3}(\Delta, \nu)$ by
\begin{displaymath}
e_{3} (\Delta, \nu) = \max \{ |E(G)| : \Delta(G)\le\Delta, \nu(G)\le\nu \}\le (\Delta+1)\nu,
\end{displaymath}
\begin{displaymath}
E_{3}(\Delta, \nu) = \{ G : \Delta(G)\le\Delta, \nu(G)\le\nu, |E(G)| = e_{3} (\Delta, \nu) \}.
\end{displaymath}

The problem of finding the maximum number of edges in a graph with fixed $\Delta$, $\nu$, and vertex number $n$ has been studied in \cite{AH} and \cite{CH}. The following theorem is proven in \cite{BK} and can also be inferred from a more generalized result provided in \cite{AH} and \cite{CH}:

$$ e_3(\Delta, \nu) = \Delta\nu+\left\lfloor\frac{\nu}{\left\lceil\frac{\Delta}{2}\right\rceil}\right\rfloor \left\lfloor\frac{\Delta}{2}\right\rfloor $$

Note also that:

\begin{equation} \label{matching bar 1}
e_3(\Delta,\nu) \leq \Delta \nu+\frac{\nu}{\lceil\frac{\Delta}{2}\rceil} \left\lfloor\frac{\Delta}{2}\right\rfloor = \left(\Delta + \frac{\lfloor{\frac{\Delta}{2}}\rfloor}{\lceil{\frac{\Delta}{2}}\rceil}\right)\nu,
\end{equation}

with the inequality becoming an equality if and only if $\lceil\frac{\Delta}{2}\rceil$ divides $\nu$.

Our main aim in this section is to characterize $E_{3}(\Delta,\nu)$. Let $G$ be a graph with $\Delta(G) \le \Delta$ and $\nu(G)\le \nu$.

\begin{rem}\label{delta1-nu1}
The set $E_{3}(1,\nu)$ is trivial when $\Delta=1$ or $\nu=1$.
\begin{itemize}
\item $\Delta=1$: if $G \in E_{3}(1,\nu)$ then $G$ consists of $\nu$ components where each component is $K_2$.
\item $\nu=1$, $\Delta=2$: $E_{3}(2,1) = \{ K_3 \}$.
\item $\nu=1$, $\Delta = 3$: $E_{3}(3,1) = \{ K_3, K_{1,3} \}$.
\item $\nu=1$, $\Delta > 3$: $E_{3}(\Delta,1) = \{ K_{1,\Delta} \}$.
\end{itemize}
\end{rem}

We next consider cases where $\Delta\geq 2$ and $\nu \geq 2$.

For a given $\Delta \geq 2$, define a graph $J_{\Delta}$ as follows:
\begin{itemize}
\item If $\Delta$ is even, $J_{\Delta}=K_{\Delta+1}$.\\
Note that $\Delta(J_{\Delta})=\Delta$, $\nu(J_{\Delta})=\frac{\Delta}{2}$, and $|E(J_{\Delta})|=\frac{(\Delta+1)\Delta}{2}$.
\item If $\Delta$ is odd, let $\Delta=2j-1$ for some $j \geq 2$. Starting with $K_{2j}$, remove a maximum matching. Connect $2j-1$ of the vertices to a new vertex to obtain the graph $J_{\Delta}$.\\
Note that $\Delta(J_{\Delta})=\Delta$, $\nu(J_{\Delta})=j=\left\lceil\frac{\Delta}{2}\right\rceil$, and, by looking at the vertex degrees, 
\begin{eqnarray*}
|E(J_{\Delta})|&=&\frac{2j\Delta+(\Delta-1)}{2}=(\Delta+\frac{(\frac{\Delta-1}{2})}{j})j=(\Delta+\frac{\lfloor{\frac{\Delta}{2}}\rfloor}{\lceil{\frac{\Delta}{2}}\rceil})\lceil{\frac{\Delta}{2}}\rceil.
\end{eqnarray*}
\end{itemize}
%\begin{minipage}{7cm}
\begin{center}
\includegraphics[scale=1.3]{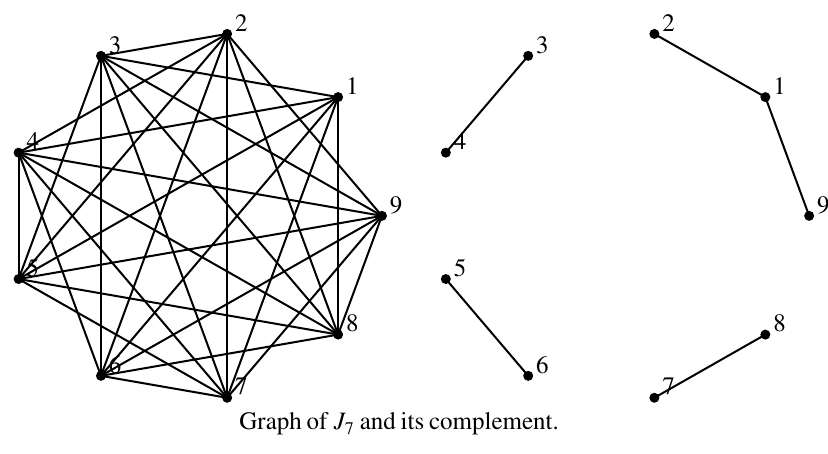}
\end{center}
%\end{minipage}
\begin{prop} \label{unique c}
Let $\Delta \geq 2$ and let $G$ be a simple graph such that $\Delta(G)=\Delta$, $\nu(G)=\left\lceil\frac{\Delta}{2}\right\rceil$, and $|E(G)| = (\Delta+\frac{\lfloor{\frac{\Delta}{2}}\rfloor}{\lceil{\frac{\Delta}{2}}\rceil})\lceil{\frac{\Delta}{2}}\rceil$. Then\\
(a) $\nu(G\setminus x)=\nu(G)$ for all $x \in V(G)$,\\
(b) $G$ is connected,\\
(c) $G \simeq J_{\Delta}$.
\end{prop}
{\sl Proof:} Let $\Delta \geq 2$. Let $G$ be a graph satisfying the conditions of the proposition.\newline
{\sl Proof of (a)}: If the statement (a) is false then there exists a vertex $x \in V(G)$ such that $\nu(G\setminus x)<\nu(G)$. As at most one edge can cover $x$ in any maximum matching, we have $\nu(G\setminus x)=\nu(G)-1$. Therefore,
\begin{displaymath}
\begin{array}{rcll}
|E(G)| &\leq & |E(G\setminus x)|+\Delta(G\setminus x)\\
& \leq & (\Delta(G\setminus x)+\frac{\lfloor{\frac{\Delta(G\setminus x)}{2}}\rfloor}{\lceil{\frac{\Delta(G\setminus x)}{2}}\rceil})\nu(G\setminus x)+\Delta(G\setminus x) & \textrm{ by equation (\ref{matching bar 1})}\\
& \leq & (\Delta+\frac{\lfloor{\frac{\Delta}{2}}\rfloor}{\lceil{\frac{\Delta}{2}}\rceil})\nu(G\setminus x)+\Delta & \textrm{ since } \Delta\geq \Delta(G\setminus x) \\
&=& (\Delta+\frac{\lfloor{\frac{\Delta}{2}}\rfloor}{\lceil{\frac{\Delta}{2}}\rceil})(\nu(G)-1)+\Delta\\
&=& (\Delta+\frac{\lfloor{\frac{\Delta}{2}}\rfloor}{\lceil{\frac{\Delta}{2}}\rceil})\nu(G)-\frac{\lfloor{\frac{\Delta}{2}}\rfloor}{\lceil{\frac{\Delta}{2}}\rceil} \\
&=& |E(G)| - \frac{\lfloor{\frac{\Delta}{2}}\rfloor}{\lceil{\frac{\Delta}{2}}\rceil} \\
&<& |E(G)|,
\end{array}
\end{displaymath}
which is a contradiction. Hence statement (a) holds.

{\sl Proof of (b)}: On the contrary assume that $G$ is not connected. Since $G$ has no isolated vertices, all components of $G$ are nontrivial, i.e., have at least an edge. Let $C_{1}$ be a component of $G$. Then $1\leq \nu(C_{1})<\nu(G)=\lceil{\frac{\Delta}{2}}\rceil $. By statement (a) and Gallai's Lemma (Lemma \ref{Gallai's lemma}), $C_{1}$ is a factor-critical component. Therefore, $|V(C_{1})|=2\nu(C_{1})+1$. So,
\begin{equation}\label{eq b}
 |E(C_{1})|\leq (2\nu(C_{1})+1)\frac{\Delta(C_{1})}{2}\leq (2\nu(C_{1})+1)\nu(C_{1}).
\end{equation}
The above inequality implies that
\begin{displaymath}
\frac{|E(C_{1})|}{\nu(C_{1})} \leq 2 \nu(C_{1})+1 \leq 2(\lceil \frac{\Delta}{2} \rceil-1)+1 < \Delta+ \frac{\lfloor\frac{\Delta}{2}\rfloor}{\lceil\frac{\Delta}{2}\rceil},
\end{displaymath}
since $2 (\lceil{\frac{\Delta}{2}}\rceil-\frac{\Delta}{2})-1 \leq 0< \frac{\lfloor\frac{\Delta}{2}\rfloor}{\lceil\frac{\Delta}{2}\rceil} $ when $\Delta \geq 2$.\\
So there is a component $C_{2}$ of $G$ such that $\frac{|E(C_{2})|}{\nu(C_{2})}>\Delta+ \frac{\lfloor\frac{\Delta}{2}\rfloor}{\lceil\frac{\Delta}{2}\rceil}$ as $\frac{|E(G)|}{\nu(G)}=\Delta+ \frac{\lfloor\frac{\Delta}{2}\rfloor}{\lceil\frac{\Delta}{2}\rceil}$. But equation (\ref{matching bar 1}) demands that $\frac{|E(C_{2})|}{\nu(C_{2})}\leq \Delta+ \frac{\lfloor\frac{\Delta}{2}\rfloor}{\lceil\frac{\Delta}{2}\rceil}$. The contradiction implies that statement (b) holds.

{\sl Proof of (c)}: Since statements (a) and (b) hold for $G$, $G$ is factor-critical by Gallai's Lemma (Lemma \ref{Gallai's lemma}). As $\nu(G)=\lceil{\frac{\Delta}{2}}\rceil$, we have $|V(G)|=2(\lceil{\frac{\Delta}{2}}\rceil)+1$. We have the following two cases.
\begin{enumerate}
\item [(i)] If $\Delta$ is even then $G$ is a connected graph with $2(\lceil{\frac{\Delta}{2}}\rceil)+1=\Delta+1$ vertices and $|E(G)|=(\Delta+ \frac{\lfloor\frac{\Delta}{2}\rfloor}{\lceil\frac{\Delta}{2}\rceil})\lceil\frac{\Delta}{2}\rceil= \frac{(\Delta+1)\Delta}{2}$. Therefore, $deg_G(x)=\Delta$ for all $x \in V(G)$. Hence $G$ is $K_{\Delta+1}$, the complete graph on $\Delta+1$ vertices. So $G\simeq J_{\Delta}$.
\item [(ii)] If $\Delta$ is odd, let $\Delta=2j-1$ for some $j\geq 2$. Then
 $\nu({G})=\left\lceil\frac{\Delta}{2}\right\rceil=j$, 
$|V(G)|=2\nu(G)+1=2j+1$ and 
$|E({G})|=(\Delta+ \frac{\lfloor\frac{\Delta}{2}\rfloor}{\lceil\frac{\Delta}{2}\rceil})\lceil\frac{\Delta}{2}\rceil=(2j-1)j+j-1.$
So \[\sum_{x \in V(G)}deg_G(x)=2j(2j-1)+2j-2.\] Therefore there is a unique vertex $v \in V({G})$ of degree $2j-2$. Hence there is a vertex $u \in V({G})$ which is not a neighbor of $v$. Consequently ${G} \setminus u$ is a regular graph of degree $2j-2$ on $2j$ vertices and hence its complement is a regular graph of degree one, namely, a matching of a complete graph on $2j$ vertices. This establishes ${G}\simeq J_{\Delta}$.\qed
\end{enumerate}

\begin{prop}\label{unique_graphs_a}
Let $\Delta, \nu \in \mathbb{Z}^{+}$ with $\Delta, \nu \geq 2$ and let $G \in E_{3}(\Delta, \nu)$. If $\lceil{\frac{\Delta}{2}}\rceil$ divides $\nu$, then $\nu(G\setminus x)=\nu$ for all $x \in V(G)$.
\end{prop}
{\sl Proof:} Let $\Delta, \nu \geq 2$ be integers and let $G \in E_{3}(\Delta, \nu)$. Suppose there exists a vertex $v \in V(G)$ such that $\nu(G\setminus v)<\nu(G)$. Then $\nu(G\setminus v)=\nu-1$ and $\Delta(G \setminus v) \leq \Delta$. This implies
\begin{eqnarray*}
|E(G)|&\leq& deg_G(v)+|E(G\setminus v)| \leq \Delta + e_{3}(\Delta(G\setminus v), \nu-1).
\end{eqnarray*}
When $\nu$ is fixed, $e_{3}(\Delta, \nu)$ is a nondecreasing function of $\Delta$. Also, $\Delta(G \setminus v) \leq \Delta$, so we now have
\begin{eqnarray*}
|E(G)|&\leq& \Delta + e_{3}(\Delta, \nu-1) \leq \Delta + \Delta(\nu-1) + \left\lfloor \frac{\nu-1}{\lceil\frac{\Delta}{2}\rceil}\right\rfloor{\left\lfloor \frac{\Delta}{2}\right\rfloor} \\
&\leq& \Delta\nu + \frac{(\nu-1)}{\lceil\frac{\Delta}{2}\rceil}{\left\lfloor \frac{\Delta}{2}\right\rfloor} = \Delta\nu + \frac{\nu}{\lceil\frac{\Delta}{2}\rceil}{\left\lfloor \frac{\Delta}{2}\right\rfloor} - \frac{{\left\lfloor\frac{\Delta}{2}\right\rfloor}}{\lceil\frac{\Delta}{2}\rceil} =|E(G)| - \frac{{\left\lfloor\frac{\Delta}{2}\right\rfloor}}{\lceil\frac{\Delta}{2}\rceil} \\
&<&  |E(G)|,
\end{eqnarray*}
which is a contradiction. Therefore no such $v$ can exist. \qed
\begin{prop}\label{unique_graphs_bc}
Let $\Delta, \nu \in \mathbb{Z}^{+}$ with $\Delta, \nu \geq 2$ and let $G \in E_{3}(\Delta, \nu)$. If $\lceil{\frac{\Delta}{2}}\rceil$ divides $\nu$ and $C$ is a component of $G$, then $\nu(C) = \lceil{\frac{\Delta}{2}}\rceil$.
\end{prop}
{\sl Proof}: Let $\Delta, \nu \geq 2$ be integers and let $G \in E_{3}(\Delta, \nu)$. Since $\lceil{\frac{\Delta}{2}}\rceil$ divides $\nu$, $e_{3}(\Delta,\nu)=(\Delta+\frac{\lfloor{\frac{\Delta}{2}}\rfloor}{\lceil{\frac{\Delta}{2}}\rceil})\nu$. Let $C$ be a component of $G$ such that $\nu(C) \neq \lceil{\frac{\Delta}{2}}\rceil$. Proposition \ref{unique_graphs_a} and Gallai's Lemma imply that $C$ is factor-critical. Hence $|V(C)|=2\nu(C)+1$. Also, for a fixed $\nu$, $e_{3}(\Delta,\nu)$ is a nondecreasing function in $\Delta$. We have the following two cases:
\begin{enumerate}
\item [(i)] If $\nu(C) > \lceil{\frac{\Delta}{2}}\rceil$, then $|E(C)| \leq \left\lfloor\frac{(2\nu(C)+1)\Delta}{2}\right\rfloor = \nu(C)\Delta +\left\lfloor\frac{\Delta}{2}\right\rfloor$ and therefore
\begin{displaymath}
\frac{|E(C)|}{\nu(C)} \leq \Delta +\frac{\lfloor\frac{\Delta}{2}\rfloor}{\nu(C)}< \Delta+\frac{\lfloor\frac{\Delta}{2}\rfloor}{\lceil\frac{\Delta}{2}\rceil}=\frac{|E(G)|}{\nu}.
\end{displaymath}
\item [(ii)] If $\nu(C) < \lceil{\frac{\Delta}{2}}\rceil$, then $|E(C)| \leq \frac{(2\nu(C)+1)2\nu(C)}{2}=(2\nu(C)+1)\nu(C)$ and therefore
\begin{eqnarray*}
\frac{|E(C)|}{\nu(C)} &\leq& 2\nu(C)+1 \leq 2(\lceil{\frac{\Delta}{2}}\rceil-1)+1\\ 
&=& 2\lceil{\frac{\Delta}{2}}\rceil - 1 < \Delta+\frac{\lfloor\frac{\Delta}{2}\rfloor}{\lceil\frac{\Delta}{2}\rceil} = \frac{|E(G)}{\nu}.
\end{eqnarray*}
\end{enumerate}

Both cases imply that there is another component $C_2$ of $G$ such that $\frac{|E(C_2)|}{\nu(C_2)} > \Delta+\frac{\lfloor\frac{\Delta}{2}\rfloor}{\lceil\frac{\Delta}{2}\rceil}$. This means $|E(C_2)| > (\Delta+\frac{\lfloor\frac{\Delta}{2}\rfloor}{\lceil\frac{\Delta}{2}\rceil})\nu(C_2) = e_{3}(\Delta, \nu(C_2))$, which is a contradiction. Therefore, all components C of $G$ satisfy $\nu(C)={\lceil\frac{\Delta}{2}\rceil}$.
\qed
\begin{prop}\label{unique_graphs_d}
Let $\Delta, \nu \in \mathbb{Z}^{+}$ with $\Delta, \nu \geq 2$ and let $G \in E_{3}(\Delta, \nu)$. If $\lceil{\frac{\Delta}{2}}\rceil$ divides $\nu$, then every component of $G$ is isomorphic to $J_{\Delta}$.
\end{prop}
{\sl Proof}: Let $\Delta, \nu \geq 2$ be integers and let $G \in E_{3}(\Delta, \nu)$. Let $C$ be a component of $G$. By Proposition \ref{unique_graphs_bc}, $\nu(C)={\lceil\frac{\Delta}{2}\rceil}$. Also, the proof of Proposition \ref{unique_graphs_bc} implies that $\frac{|E(C)|}{\nu(C)} = \frac{|E(G)}{\nu}$, so $|E(C)|=(\Delta+\frac{\lfloor\frac{\Delta}{2}\rfloor}{\lceil\frac{\Delta}{2}\rceil}){\lceil\frac{\Delta}{2}\rceil}$. By Proposition \ref{unique c}, $C$ must be isomorphic to $J_{\Delta}$. \qed
\begin{prop}\label{inverse_unique_graph}
Let $\Delta, \nu \in \mathbb{Z}^{+}$ and let $G_1 \in E_{3}(\Delta,\nu)$. If $\Delta, \nu \geq 2$ and ${\lceil{\frac{\Delta}{2}}\rceil}$ doesn't divide $\nu$, then there exists a simple graph $G_2$ such that $G_2 \in E_{3}(\Delta,\nu)$ and $G_2$ is \underline{not} isomorphic to $G_1$. 
\end{prop}

{\sl Proof:} Let $\Delta, \nu \in \mathbb{Z}^{+}$ such that $\Delta, \nu \geq 2$ and ${\lceil{\frac{\Delta}{2}}\rceil}$ doesn't divide $\nu$. Let $G_1 \in E_{3}(\Delta,\nu)$. We use the method given in \cite{BK} to construct a graph $G_2 \in E_{3}(\Delta,\nu)$ that is not isomorphic to $G_1$. 

Let $s:=\left\lfloor{\frac{\nu(G_1)}{{\lceil{\frac{\Delta(G_1)}{2}}\rceil}}}\right\rfloor$ and $t:=\nu(G_1)-({\lceil{\frac{\Delta(G_1)}{2}}\rceil})s$. Let $G$ be the graph with $s+t$ components where $s$ components are isomorphic to $J_{\Delta}$ and $t$ components are isomorphic to $K_{1,\Delta}$. If $G_1$ is \underline{not} isomorphic to $G$, then take $G_2=G$. Otherwise, assume $G_1$ is isomorphic to $G$. 

Since ${\lceil{\frac{\Delta}{2}}\rceil}$ doesn't divide $\nu$, $t \geq 1$. If $t \geq 2$, let $G_2$ be the graph with $s+t-1$ components where $s$ components are isomorphic to $J_{\Delta}$, $t-2$ components are isomorphic to $K_{1,\Delta}$, and one component is isomorphic to $K_{2,\Delta}$. It is clear by construction that $G_2 \in E_{3}(\Delta,\nu)$ and $G_2$ is not isomorphic to $G_1$. 

Now suppose $t=1$. Since $\nu \geq 2$ and $\nu(K_{1,\Delta})=1$, $G_1$ must have another component, i.e. $s \geq 1$ and $G_1$ has a component isomorphic to $J_{\Delta}$. Since

\begin{displaymath}
\Delta+|E(J_{\Delta})|=\Delta+\Delta\left\lceil{\frac{\Delta}{2}}\right\rceil+\left\lfloor{\frac{\Delta}{2}}\right\rfloor \leq \left\lfloor{\frac{(2(\lceil\frac{\Delta}{2}\rceil +1)+1)\Delta}{2}}\right\rfloor,
\end{displaymath}

we can merge the two components to form a factor critical component $C$ such that $|V(C)|=2(\lceil\frac{\Delta}{2}\rceil +1)+1$, $|E(C)|=\Delta+|E(J_{\Delta})|$, $\Delta(C)=\Delta$, and $\nu(C)=\nu(J_{\Delta})+1$. Then take $G_2$ to be the graph with $s$ components where $s-1$ components are isomorphic to $J_{\Delta}$ and one component is isomorphic to $C$. It is clear by construction that $G_2 \in E_{3}(\Delta,\nu)$ and $G_2$ is not isomorphic to $G_1$. \qed

\begin{thm}\label{uniqueness_thm}
Let $\Delta, \nu \in \mathbb{Z}^{+}$ such that $\Delta, \nu \geq 2$ and  and let $G \in E_{3}(\Delta,\nu)$ such that $G$ has no isolated vertices. $G$ is a unique graph up to isomorphism if and only if ${\lceil{\frac{\Delta(G)}{2}}\rceil}$ divides $\nu(G)$.
\end{thm}
{\sl Proof:} This follows directly from equation (\ref{matching bar 1}), Proposition \ref{unique_graphs_d}, and Proposition \ref{inverse_unique_graph}.\qed

Now Theorem \ref{delta-nu2} is an easy consequence of Remark \ref{delta1-nu1} and Theorem \ref{uniqueness_thm}.
%****************************************************************************************
%********************************************************************************************
\section{Acknowledgements}
We would like to thank Prof. \'Akos Seress for his valuable comments.    


\begin{thebibliography}{99}
\bibitem[1]{AF} J. Akiyama and P. Frankl, On the size of graphs with complete-factors, J. of Graph Theory, \textbf{9} (1985) 197-201.
\bibitem[2]{AH} H.L. Abbott, D. Hanson and N. Sauer, Intersection theorems for systems of sets, J. Combinatorial Theory (A) \textbf{12} (1972), 381-389.
\bibitem[3]{CH} V. Chv\'atal and D. Hanson, Degrees and Matchings, J. Combinatorial Theory (B) \textbf{20} (1976), 128-138. 
\bibitem[4]{BK} N. Balachndran and N. Khare, Graphs with restricted valency and matching number, Discrete Mathematics, \textbf{309} (2009), 4176-4180.
\bibitem[5]{EG} P. Erd\H{o}s and T. Gallai, On the mininal number of vertices representing the edges of a graph, Publ. Math. Inst. Hungar. Acad. Sci, \textbf{6} (1961), 181-203.
\bibitem[6]{ER} P. Erd\H{o}s and R. Rado, Intersection theorems for systems of sets, J. London Math.Soc., \textbf{35} (1960), 85-90.
\bibitem[7]{LP} L. Lovasz and M. D. Plummer, Matching Theory, North Holland, 1986.
\bibitem[8]{DW} D. B. West, Introduction to Graph Theory, Prentice Hall, 1996.
\end{thebibliography}
\end{document}